\documentclass[12pt,draftcls,onecolumn]{IEEEtran}

\usepackage{amssymb,latexsym,color,amsmath,pifont,epsfig,graphicx}
\usepackage{amsthm}
\usepackage{setspace,cite}
\usepackage{subfig}  
\usepackage{url}

\usepackage{tikz}
\usepackage{pgfplots,pgfplotstable}
\graphicspath{{}} 

\newcommand{\mre}{\mathrm{e}}
\newcommand{\mrd}{\mathrm{d}}
\newcommand{\ds}{\displaystyle}

\newcommand{\DefinedAs}[0]{\mathrel{\mathop:}=}

\newcommand{\R}{\mathbb{R}}

\newtheorem{lemma}{Lemma}
\newtheorem{proposition}{Proposition}

\IEEEoverridecommandlockouts                              

\begin{document}

\title{\LARGE \bf Performance of leader-follower multi-agent systems in directed networks}

\author{Fu Lin
\thanks{F.\ Lin is with the Systems Department, United Technologies Research Center, 411 Silver Ln, East Hartford, CT 06118, USA (e-mail: linf@utrc.utc.com).}
\thanks{This material is based upon work supported by the U.S. Department of Energy, Office of Science, Office of Advanced Scientific Computing Research, Applied Mathematics program under contract number DE-AC02-06CH11357.}}

\maketitle

    \begin{abstract}
We consider leader-follower multi-agent systems in which the leader executes the desired trajectory and the followers implement the consensus algorithm subject to stochastic disturbances. The performance of the leader-follower systems is quantified by using the steady-state variance of the deviation of the followers from the leader. We study the asymptotic scaling of the variance in directed lattices in one, two, and three dimensions. We show that in 1D and 2D the variance of the followers' deviation increases to infinity as one moves away from the leader, while in 3D it remains bounded. We prove that the variance scales as a square-root function in 1D and a logarithmic function in 2D lattices.
    \end{abstract}
    

\section{Introduction}

A leader-follower multi-agent system consists of a leader, who provides the desired trajectory of the multi-agent system, and a set of followers, who update their states using local relative feedback. This control strategy has a variety of applications including formation of unmanned air vehicles, control of rigid robotic bodies, and distributed estimation in sensor networks~\cite{barhes06,barhes07,barhes08,rahjimesege09,youscaleo10,linfarjovCDC11,farlinjovCDC11,clapoo11,bamjovmitpat12,linfarjovTAC12platoons,linfarjov12,clabuspoo14}.

A fundamental question concerning the performance of the leader-follower strategy is how well the followers are able to keep track the trajectory of the leader when they are subject to stochastic disturbances. In large networks, the asymptotic scaling of the variance of followers' deviation from the desired trajectory is determined by the network architecture. In this paper, we focus on directed lattices in one, two, and three dimensions. We show that as one moves away from the leader, the variance of the followers increases unboundedly in 1D and 2D, whereas in 3D the variance of the followers is bounded above by a constant that is independent of the number of followers. These results resemble the performance limitation of distributed consensus in undirected tori~\cite{bamjovmitpat12}. For {\em directed networks\/}, our results for the asymptotic scaling of the performance appear to be among the first in the literature. 

Our contributions are twofold. First, we obtain analytical expressions for the steady-state variance of the deviation of the followers from the leader. These expressions allow us to study the distribution of variance in leader-follower multi-agent systems with directed lattices as the controller architecture. Second, we characterize the asymptotic scaling trends of the variance of the followers in 1D, 2D, and 3D directed lattices. We show that in 1D and 2D the variance of the followers scales asymptotically as a square-root function and a logarithmic function, respectively, and in 3D the variance remains bounded regardless of the network size.

This paper is organized as follows. In Section~\ref{sec.var}, we present our main results for the performance of leader-follower multi-agent systems on directed lattices. We also discuss connection between our results and random walks on undirected lattices. In Section~\ref{sec.proof} we provide the proofs and in Section~\ref{sec.concl} we summarize our findings.

\section{Leader-follower multi-agent systems on directed lattices}
\label{sec.var}

We consider the performance of leader-follower multi-agent systems on directed lattices. By exploiting the lower triangular Toeplitz structure of the modified Laplacian matrices, we obtain analytical expressions for the variance of followers and establish its asymptotic scaling trends in large networks. 	 


     \begin{figure}
       \centering
         \includegraphics[width=0.75\textwidth]{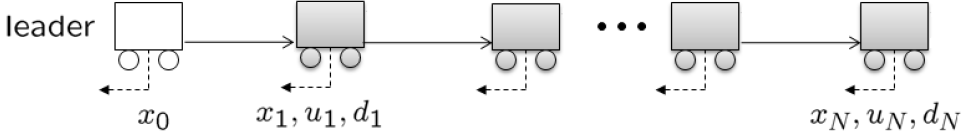}
       \caption{The leader-follower system on 1D lattice. The first follower has access to the state of the leader indexed by $0$.}
       \label{fig.1D}
     \end{figure}

\subsection{1D lattice}

Consider a set of $N$ agents on a line whose dynamics are modeled by the single integrators
\[
  \dot{x}_n(t) \,=\, u_n(t) \,+\, d_n(t),
  \quad n \,=\, 1, \ldots, N,
\]
where $x_n(t)$ denotes the deviation of the $n$th vehicle from its desired trajectory, $u_n(t)$ is the control input, and $d_n(t)$ is a zero-mean, unit-variance stochastic disturbance. A virtual leader, indexed by $0$, is assumed to execute the desired trajectory at all times. Thus, its deviation from the desired trajectory is zero,
$
  x_0(t) \,\equiv \, 0,
$ 
and $\dot{x}_0(t) = 0$. The followers implement the consensus algorithm. Namely, each follower updates its state information using the relative differences between itself and the agent ahead (see Fig.~\ref{fig.1D}):
\[ 
  \dot{x}_n(t) \,=\, - \, (x_n(t) \,-\, x_{n-1}(t)) \,+\, d_n(t),
  \quad n \,=\, 1, \ldots, N.
\]
We assume that the first follower has access to the state of the leader. Since $x_0(t) \equiv 0$, it follows that
\[
  \dot{x}_1(t) \,=\, - x_1(t) \,+\, d_1(t).
\] 
By stacking the states of all followers into a vector, $x(t) = [ \, x_1(t) \, \cdots \, x_N (t) \, ]^T \in \R^N$, the state-space representation of the leader-follower system is given by
\begin{equation} 
 \label{eq.1dss} 
  \dot{x}(t) \,=\, - L x(t) \,+\, d(t),
\end{equation}
where $L \in \R^{N \times N}$ is the modified Laplacian matrix of the 1D lattice. In particular, $L$ is lower triangular Toeplitz with $1$ on the main diagonal, $-1$ on the first subdiagonal, and zero everywhere else:
\begin{equation}
  \label{eq.L1d}
        L
        \, = \,
        \left[
        \begin{array}{rrrr}
        1 & 0 & \cdots  & 0  \\
        -1 & 1 & \ddots & \vdots  \\
        0 & \ddots & \ddots & 0 \\
        0 & \cdots & -1 & 1 
        \end{array}
        \right].
\end{equation}

When the disturbance, $d(t) = [ \, d_1(t) \, \cdots \, d_N(t) ]^T \in \R^N$, is absent, the deviation of the followers asymptotically converges to zero. In other words, the followers converge to the desired trajectory, that is, the trajectory of the leader. In the presence of the disturbance, however, the followers converge to the desired state in the mean value. The steady-state variance of the followers can be used to quantify the deviation from the desired state:
\[
  V_n \,\DefinedAs \, \lim_{t \to \infty} E\{x_n^2(t)\},
  \quad 
  n \,=\, 1,\ldots,N,
\]
where $E\{ \cdot \}$ denotes the expectation operator.

We are interested in the scaling trend of the variance distribution as one moves away from the leader. Intuitively, the followers who are farther away from the leader have larger steady-state variance. It turns out that the variance of the followers increases as a square-root function of the number of followers. This result is detailed in Lemma~\ref{lem.1d}.

\begin{lemma}
\label{lem.1d}
The steady-state variance of the $n$th follower in the 1D lattice~\eqref{eq.1dss} is given by
\begin{equation}
  \label{eq.pn}
       V_{n}
        \,=\,
        \sum_{i=1}^n \frac{(2i-2)!}{2\cdot 2^{2i-2}((i-1)!)^2}
        \,=\,
        \frac{n \, (2n)!}{2^{2n} \, n! \, n!},
        \;\; n = 1, \ldots, N.
\end{equation}
The total variance normalized by the number of followers is 
    \[
        \Pi_N
        \, \DefinedAs \,
        \frac{1}{N} \, \sum_{n=1}^N V_n
        \,=\,
        \frac{(2N+1)!}{3 \cdot 2^{2N}  N! N!}.
    \] 
Furthermore,
\begin{equation}
   \nonumber
    \lim_{n \to \infty}
    \frac{V_{n}}{\sqrt{n}}
    \, = \,
    \sqrt{\frac{1}{\pi}},
    ~~~
    \ds{\lim_{N  \to \infty}}
    \dfrac{\Pi_N}{\sqrt{N}}
    \, = \,
    \dfrac{2}{3 \sqrt{\pi}}.
\end{equation}
\end{lemma}
The proof can be found in Section~\ref{app.1d}.

\begin{figure}
       \centering
         \includegraphics[width=0.4\textwidth]{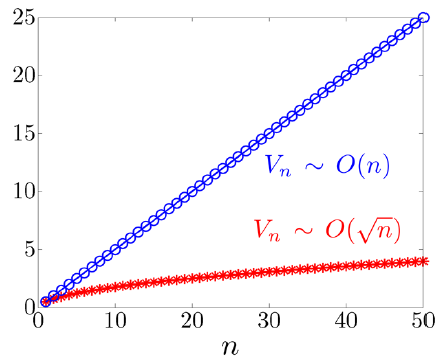}
       \caption{The variance of followers in 1D scales as square-root and linear functions for directed ($\textcolor{red}{*}$) and undirected ($\textcolor{blue}{\circ}$) lattices, respectively.}
       \label{fig.Vn1D}
     \end{figure}

To put Lemma~\ref{lem.1d}  in context, recall that the variance of the undirected 1D lattice scales as a linear function of $n$; see e.g.,~\cite{linfarjovTAC12platoons,linfarjov12}. This implies that the control architecture with directed networks outperforms the undirected counterpart in 1D lattices; see Fig.~\ref{fig.Vn1D}.

\subsection{2D lattice}

We next consider the leader-follower system that consists of a virtual leader and $N \times N$ followers in the formation of a 2D lattice. A follower at the $n$th row and the $m$th column of the 2D lattice, indexed by $(n,m)$, updates its state using the relative differences between itself and its two neighbors:
\[
  \dot{x}_{n,m} 
  \,=\, 
  - \, (x_{n,m} - x_{n,m-1}) 
  \,-\, (x_{n,m} - x_{n-1,m}) 
  \,+\, d_{n,m},
\]
for $n,m=1,\ldots,N$. Here, we drop the dependence on time in order to ease the notation. Recall that in the 1D case, the first follower in the formation is assumed to have access to the state of the leader. Similarly, we assume that the followers on the boundary of the 2D formation have direct access to the state of the leader. In particular, the followers on the first column and the first row implement the following closed-loop dynamics
\[
  \begin{array}{l}
  \dot{x}_{n,1} \,=\, - \, (x_{n,1} - x_{n-1,1}) \,-\, (x_{n,1} - x_{n,0}) \,+\, d_{n,1}, \\[0.1cm]
  \dot{x}_{1,m} \,=\, - \, (x_{1,m} - x_{0,m})  \,-\, (x_{1,m} - x_{1,m-1}) \,+\, d_{1,m},
  \end{array}
\]
where $x_{n,0} = x_{0,m} = x_0 \equiv 0$. 

Let ${\bf x} = [{\bf x}_1^T \cdots {\bf x}_N^T]^T \in \mathbb{R}^{N^2}$ be the state of followers where  ${\bf x}_n \,=\, [\, x_{n,1} \,\cdots \, x_{n,N} \, ]^T \in \mathbb{R}^{N}$ denotes the state of followers on the $n$th row of the lattice. Then the state-space representation of the leader-follower system is given by
\begin{equation}
   \label{eq.2dss} 
        \dot{{\bf x}}
        \, = \,
        - \, L_2 {\bf x}
        \, + \,
        d,
\end{equation}
where the modified Laplacian matrix $L_2 \in \R^{N^2 \times N^2}$ is lower triangular block Toeplitz:
    \begin{equation}
    \nonumber
        L_2
        \, = \,
        \left[
        \begin{array}{rrrr}
        K_2 & 0 & \cdots  & 0  \\
        -I & K_2 & \ddots & \vdots  \\
        0 & \ddots & \ddots & 0 \\
        0 & \cdots & -I & K_2
        \end{array}
        \right],
    \end{equation}
where $I$ is the identity matrix and $K_2 \in \R^{N \times N}$ is lower triangular Toeplitz with $2$ on its main diagonal, $-1$ on the first subdiagonal, and zero everywhere else:
    \[
        K_2
        \, = \,
        \left[
        \begin{array}{rrrr}
        2 & 0 & \cdots  & 0  \\
        -1 & 2 & \ddots & \vdots  \\
        0 & \ddots & \ddots & 0 \\
        0 & \cdots & -1 & 2
        \end{array}
        \right].
    \]

In what follows, we derive the analytical expression for the variance of each follower. The steady-state covariance matrix of the leader-follower system is given by the solution of the Lyapunov equation
\[
(- L_2) P \,+\, P (- L_2)^T
\,=\, I.
\]
Alternatively, the covariance matrix can be expressed by the integral form
    \[
        P \,=\,
        \int_0^\infty
        \mre ^{ - L_2 t }
        \,
        \mre ^{ - L_2^T t }
        \, \mrd t
        \, \in \,
        \mathbb{R}^{N^2 \times N^2}.
    \]
Let $P_n \in \mathbb{R}^{N \times N}$ be the $n$th diagonal block of $P$, and let $(P_{n})_{m}$ be the $m$th diagonal element of $P_n$  for $n, m = 1,\ldots,N$. We have the following result.
\begin{lemma}
\label{lem.2d}
For the leader-follower system in the 2D lattice~\eqref{eq.2dss}, the steady-state variance of the follower at the $n$th row and $m$th column is given by
\begin{equation}
  \label{eq.pnm}
        (P_{n})_{m}
        \, = \,
        \sum_{i  =  1}^{n}
        \sum_{j  =  1}^{m}
        \frac{ (2i + 2j - 4)! }{4 \cdot 4^{2i+2j-4} ( (i-1)! (j-1)! )^2}
\end{equation}
for $n,m=1,\ldots,N$.
\end{lemma}
The proof of Lemma~\ref{lem.2d} can be found in Section~\ref{app.2d}. Note the resemblance of the expression in the double summation~\eqref{eq.pnm} and in the single summation~\eqref{eq.pn} for the 2D and 1D lattices, respectively.

\begin{figure}
       \centering
         \includegraphics[width=0.35\textwidth]{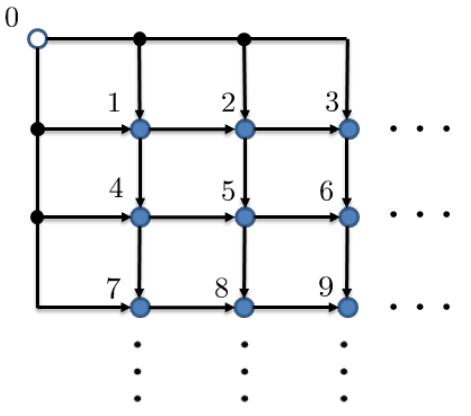}
	\quad
         \includegraphics[width=0.4\textwidth]{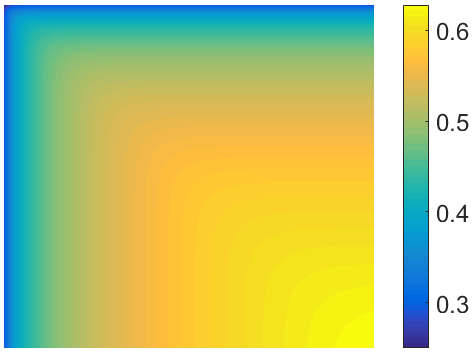}
       \caption{The leader-follower system in 2D lattice (left) and the variance of the followers of a $50 \times 50$ lattice (right). As one moves away from the leader index by $0$, the variance of the followers increases.}
       \label{fig.2D}
     \end{figure}

Since we are summing up positive quantity in~\eqref{eq.pnm}, we conclude that $(P_n)_m$ is monotonically increasing as both $n$ and $m$ increase; see Fig.~\ref{fig.2D}. In other words, the variance of the follower grows as one moves away from the leader. We next show that the variance of the followers on the diagonal of the lattice scales asymptotically as a {\em logarithmic\/} function.

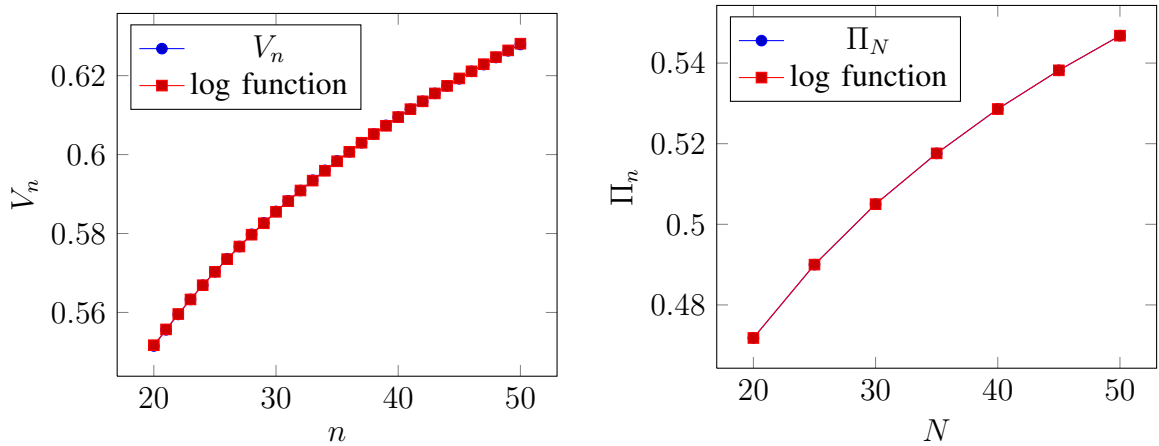
\begin{figure}
  \centering
  \begin{tikzpicture}
    \begin{axis} [width=0.45\textwidth,
      xlabel = $n$, ylabel=$V_n$,legend pos=north west
      ]
      \addplot table[x=n,y=Vn] {TwoDimVn.txt};
      \addplot table[x=n,y=logfit] {TwoDimVn.txt};
      \legend{$V_n$,log function}
    \end{axis}
  \end{tikzpicture}
  \quad
  \begin{tikzpicture}
    \begin{axis} [width=0.45\textwidth,
      xlabel = $N$, ylabel=$\Pi_n$,legend pos=north west
      ]
      \addplot table[x=Nval, y=TotalVnNorm] {TwoDimTotalVn.txt};
      \addplot table[x=Nval,y=logfit] {TwoDimTotalVn.txt};
      \legend{$\Pi_N$,log function}
    \end{axis}
  \end{tikzpicture}
  \caption{The variance of followers on the diagonal of 2D lattice and the log function $0.0834 \log(n)+0.3019$ (left); the normalized variance $\Pi_N$ and the log function $0.0819 \log(N) +   0.2263$ (right).}
  \label{fig.2Ddiag}
\end{figure}

\begin{proposition}
\label{pro.2d}
Consider the leader-follower system in the 2D lattice~\eqref{eq.2dss}. Let $V_n$ be the steady-state variance of the follower at the $n$th row and the $n$th column of the lattice for $n = 1, \ldots, N$. Then $V_n$ scales asymptotically as a logarithmic function of $n$, denoted as
$
  V_n  \, \sim \, O(\log(n)).
$
\end{proposition}
The proof can be found in Section~\ref{app.2dscale}. From Proposition~\ref{pro.2d}, it follows that the total variance of the followers on the main diagonal normalized by $N$ scales logarithmically for large $N$, that is, 
\[
\Pi_N \,\DefinedAs\,  \frac{1}{N} \sum_{n=1}^N V_n \, \sim \, O(\log(N)).
\] 
We verify Proposition~\ref{pro.2d} via numerical computation. The results are shown in Fig.~\ref{fig.2Ddiag}.

\subsection{3D lattice}

While the variance of the followers increases unboundedly with the size of lattices in 1D and 2D, it turns out that in 3D the variance of the followers is bounded by a constant that is independent of the lattice size. For undirected networks, similar results have been shown for distributed consensus~\cite{bamjovmitpat12} and distributed estimation~\cite{barhes07,barhes08}. To our best knowledge, our result for directed lattices is the first in the literature.

Consider the leader-follower system that consists of a virtual leader and  $N \times N \times N$ followers on the 3D lattice. The coordinates of the follower at the $n$th row and $m$th column of the $l$th cross section is denoted by $(n,m,l)$ for $n,m,l=1,\ldots,N$. The follower updates its state using local feedback subject to disturbance:
\begin{align}
  \nonumber
  \dot{x}_{n,m,l} 
  \,=\, 
  & - (x_{n,m,l} - x_{n-1,m,l}) \,-\, (x_{n,m,l} - x_{n,m-1,l})  \\
  \nonumber 
  & - (x_{n,m,l} - x_{n,m,l-1}) \,+\, d_{n,m,l}.
\end{align}
Similar to the 1D and 2D cases, the followers on the boundary, indexed by $(1,m,l)$, $(n,1,l)$, and $(n,m,1)$, have access to the state of the leader, that is, $x_{0,m,l} = x_{n,0,l} = x_{n,m,0} = x_0 \equiv 0$.

The state-space representation of the leader-follower system on the 3D lattice is given by
\begin{equation}
  \label{eq.3dss}
       \dot{{\bf x}}
        \, = \,
        - \, L_3 {\bf x}
        \, + \,
        d,
\end{equation}
where the modified Laplacian matrix $L \in \R^{N^3 \times N^3}$ is lower triangular block Toeplitz:
    \[
        L_3
        \, = \,
        \left[
        \begin{array}{rrrr}
        K  & 0 & \cdots  & 0  \\
        -I & K & \ddots & \vdots  \\
        0 & \ddots & \ddots & 0 \\
        0 & \cdots & -I & K
        \end{array}
        \right],
    \]
where $K \in \R^{N^2 \times N^2}$ is also lower triangular block Toeplitz:
    \[
        K
        \, = \,
        \left[
        \begin{array}{rrrr}
        K_3 & 0 & \cdots  & 0  \\
        -I & K_3 & \ddots & \vdots  \\
        0 & \ddots & \ddots & 0 \\
        0 & \cdots & -I & K_3
        \end{array}
        \right],
    \]
where $K_3 \in \R^{N \times N}$ is lower triangular Toeplitz with $3$ on the main diagonal, $-1$ on the first subdiagonal, and $0$ everywhere else:
    \[
        K_3
        \, = \,
        \left[
        \begin{array}{rrrr}
        3 & 0 & \cdots  & 0  \\
        -1 & 3 & \ddots & \vdots  \\
        0 & \ddots & \ddots & 0 \\
        0 & \cdots & -1 & 3
        \end{array}
        \right].
    \]
Similar to the 1D and 2D cases, we obtain the following expression for the steady-state variance of followers.
    \begin{lemma}
      \label{lem.3d}
Consider the leader-follower system on 3D lattice~\eqref{eq.3dss}. The steady-state variance of the follower at coordinates $(n,m,l)$ of the 3D lattice can be expressed as
\begin{equation}
  \label{eq.3d}
   ((P_n)_m)_l =  \sum_{i=1}^n   \sum_{j=1}^m   \sum_{k=1}^l 
   \frac{ 6^{6-(2i+2j+2k)} (2i+2j+2k-6)!}{6 ((i-1)!(j-1)!(k-1)!)^2}.
\end{equation}
    \end{lemma}
The proof can be found in Section~\ref{app.3d}. Note the resemblance of the expression in the triple summation~\eqref{eq.3d} and in the double summation~\eqref{eq.pnm} for the 3D and 2D lattices, respectively. 

From~\eqref{eq.3d}, we see that $((P_n)_m)_l$ is monotonically increasing as $n$, $m$, and $l$ increase. In other words, the variance of the follower grows as one moves away from the leader. Similar observations have been noted for the 1D and 2D cases. It turns out that the variance of the followers on the diagonal of the 3D lattice is bounded above by a constant independent of lattice size. 

\begin{figure}
       \centering
         \includegraphics[width=0.3\textwidth]{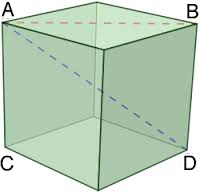}
	\qquad
         \includegraphics[width=0.4\textwidth]{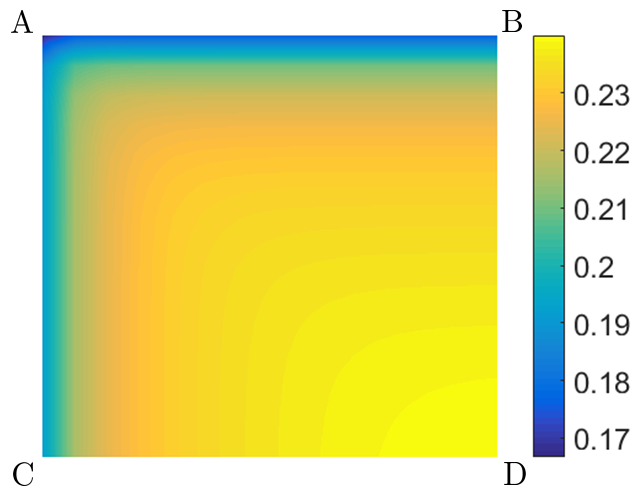}
       \caption{The diagonal slice $ABCD$ in 3D lattice (left) and the variance of the followers on the diagonal slice of the $15 \times 15 \times 15$ lattice (right). The followers on the diagonal $AD$ are denoted by coordinates $(n,n,n)$ for $n=1,\ldots,N$.}
       \label{fig.3D}
     \end{figure}

\begin{proposition}
\label{pro.3d}
Consider the leader-follower system on the 3D lattice~\eqref{eq.3dss}. Let $V_{n}$ be the steady-state variance of the follower at the coordinates $(n,n,n)$ of the 3D lattice for $n=1,\ldots,N$. Then $V_{n}$ is bounded above by a constant that is independent of network size, denoted as
$
  V_n \sim O(1).
$  
\end{proposition}
The proof is analogous to the proof of Proposition~\ref{pro.2d}; see Section~\ref{app.3dscale}. Figure~\ref{fig.3D} shows the variance of the followers on the diagonal slice of a $15 \times 15 \times 15$ lattice. The detailed numbers of the variance is provided in Table~\ref{tab.3D}. Note that the variance grows but is bounded above by constants along columns, rows, and diagonals of diagonal slice.

\begin{table}
\begin{filecontents*}{ThreeDimDiagFaceN15.txt}
    0.1667    0.1883    0.1914    0.1919    0.1920    0.1920    0.1920    0.1920    0.1920    0.1920    0.1920    0.1920    0.1920    0.1920    0.1920
    0.1759    0.2063    0.2133    0.2152    0.2158    0.2160    0.2160    0.2160    0.2160    0.2160    0.2160    0.2160    0.2160    0.2160    0.2160
    0.1767    0.2093    0.2186    0.2221    0.2235    0.2240    0.2242    0.2243    0.2243    0.2243    0.2243    0.2243    0.2243    0.2243    0.2243
    0.1768    0.2099    0.2202    0.2247    0.2269    0.2279    0.2284    0.2286    0.2287    0.2287    0.2287    0.2287    0.2287    0.2287    0.2287
    0.1768    0.2100    0.2206    0.2257    0.2285    0.2300    0.2308    0.2312    0.2314    0.2315    0.2315    0.2316    0.2316    0.2316    0.2316
    0.1768    0.2100    0.2207    0.2260    0.2292    0.2311    0.2322    0.2329    0.2332    0.2334    0.2335    0.2336    0.2336    0.2336    0.2336
    0.1768    0.2100    0.2208    0.2262    0.2294    0.2316    0.2330    0.2339    0.2344    0.2348    0.2349    0.2350    0.2351    0.2351    0.2351
    0.1768    0.2100    0.2208    0.2262    0.2295    0.2318    0.2334    0.2345    0.2352    0.2357    0.2360    0.2361    0.2362    0.2363    0.2363
    0.1768    0.2100    0.2208    0.2262    0.2296    0.2319    0.2336    0.2348    0.2357    0.2363    0.2367    0.2370    0.2371    0.2372    0.2373
    0.1768    0.2100    0.2208    0.2262    0.2296    0.2320    0.2337    0.2350    0.2360    0.2367    0.2372    0.2376    0.2378    0.2380    0.2381
    0.1768    0.2100    0.2208    0.2262    0.2296    0.2320    0.2337    0.2351    0.2361    0.2370    0.2376    0.2380    0.2383    0.2385    0.2387
    0.1768    0.2100    0.2208    0.2262    0.2296    0.2320    0.2337    0.2351    0.2362    0.2371    0.2378    0.2383    0.2387    0.2390    0.2392
    0.1768    0.2100    0.2208    0.2262    0.2296    0.2320    0.2338    0.2351    0.2363    0.2372    0.2379    0.2385    0.2389    0.2393    0.2395
    0.1768    0.2100    0.2208    0.2262    0.2296    0.2320    0.2338    0.2351    0.2363    0.2372    0.2380    0.2386    0.2391    0.2395    0.2398
    0.1768    0.2100    0.2208    0.2262    0.2296    0.2320    0.2338    0.2351    0.2363    0.2372    0.2380    0.2386    0.2392    0.2396    0.2400
\end{filecontents*}
\pgfplotstabletypeset[
	zerofill,precision=2, header=false
           ]
	{ThreeDimDiagFaceN15.txt}
\caption{Variance of the followers on the diagonal slice $ABCD$ of the $15 \times 15 \times 15$ lattice in Fig.~\ref{fig.3D}. The variance is bounded above by constants along the columns, rows, and diagonals.}
\label{tab.3D}
\end{table}

\subsection{Connections with random walks}

The connections between random walks and distributed estimation and control problems have been studied by several authors; see~\cite{barhes06,barhes07,barhes08,bamjovmitpat12,linfarjov12,fitleo13,clabuspoo14}. All existing work focuses on {\em undirected\/} networks. We next show that the asymptotic scaling for the variance of followers in {\em directed\/} lattices can be expressed as
\[
  V_n \, \sim \, \frac{1}{2D} \sum_{k=0}^{n-1} u_{2k},
\]
where $D = 1$, $2$, or $3$ is the dimension and $u_{2k}$ is the probability of a random walk of length $2k$ returning to the starting point on the undirected lattices.

Recall that for the 1D lattice, $u_{2k}$ is given by~\cite[Section 7.2]{doysne84} 
\[
        u_{2k}
        \, = \,
        \frac{1}{2^{2k}} {2k \choose k}
        \, = \,
        \frac{(2k)!}{2^{2k} \,k! \,k!}.
    \]
From~\eqref{eq.pn}, it follows that the variance of the followers in 1D can be expressed as 
\begin{equation}
  \label{eq.sum1d}
        V_{n}
        \, = \,
        \frac{1}{2}
        \sum_{k=0}^{n-1} u_{2k}.
\end{equation}
In other words, the steady-state variance of the $n$th follower can be expressed as the sum of the probability of a random walk returning to the starting point  of length $2k$ for $k=0,1,\ldots, n-1$. In 2D lattice,  $u_{2k}$ is given by~\cite[Section 7.3]{doysne84} 
    \[
        u_{2k}
        \, = \,
        \left(
        \frac{1}{2^{2k}}{2k \choose k}
        \right)^2
        \,=\, 
        \frac{1}{4^{2k}}\left( \frac{(2k)!}{k!k!} \right)^2.
    \]
From~\eqref{eq.Sk}, it follows that $S_k \, = \, (1/4) u_{2k}$. Then the summation of the positive function $f$ over the triangle ${\cal T}_n$ in the 2D lattice~\eqref{eq.triangle} can be expressed as
\begin{equation}
  \label{eq.sum2d}
  \Delta_{n} 
  \,=\, 
  \frac{1}{4} \sum_{k=0}^{n-1} u_{2k}.   
\end{equation} 
In the 3D lattice,  $u_{2k}$ is given by~\cite[Section 7.3]{doysne84}
\[
        u_{2k}
        \, = \,
         \sum_{j=0}^p \sum_{k=0}^{p-j} \frac{1}{2^{2p}} 
         \left(\frac{(2p)!}{p!p!}\right) 
         \left(\frac{p!}{3^{p} j! k! (p-j-k)!} \right)^2.
\]
It can be shown that the summation of an appropriate positive function over the triangular pyramid is given by~(see Section~\ref{app.3d})
\begin{equation}
  \label{eq.sum3d}
  T_n \,=\, \frac{1}{6} \sum_{k=0}^{n-1} u_{2k}.
\end{equation}
From~\eqref{eq.sum1d},~\eqref{eq.sum2d}, and~\eqref{eq.sum3d}, we observe that the asymptotic scaling for the variance of the followers on directed lattices can be expressed as 
\[
  V_n \, \sim \, \frac{1}{2D} \sum_{k=0}^{n-1} u_{2k},
  \quad 
  D = 1, 2, 3.
\]

\section{Proofs}
\label{sec.proof}

\subsection{Proof of Lemma~\ref{lem.1d}}
\label{app.1d}

We begin with the steady-state covariance matrix
\begin{equation}
  \label{eq.P}
  P 
  \,\DefinedAs \, \lim_{t\to\infty} E\{x(t)x^T(t)\}
  \,=\, \int_{0}^\infty \mre ^{-Lt} \mre^{-L^Tt} \mrd t.
\end{equation}
We compute the matrix exponential by using the inverse Laplace transform
$
  \mre^{-Lt}  \,=\, {\cal L}^{-1} \{ (sI \,+\, L)^{-1}  \}.
$ 
Since $L$ is a lower triangular Toeplitz matrix~(see~\eqref{eq.L1d}), it follows that  $(sI + L)^{-1}$ is also lower triangular Toeplitz 
    \[
    (sI \,+\, L)^{-1}
    \, \sim \,
    \left[
    \begin{array}{ccc}
    (s+1)^{-1} &      0     &      0       \\
    (s+1)^{-2} &(s+1)^{-1}  &      0       \\
    (s+1)^{-3} &(s+1)^{-2}  & (s+1)^{-1}   
    \end{array}
    \right].
    \]
In particular, $(s+1)^{-i}$ is the $i$th entry of the first column. By using the formula for the inverse Laplace transform
    \[ 
    {\cal L}^{-1}\{
    (s+1)^{-i}\}
    \,=\,
    \frac{t^{i-1}}{(i-1)!} \mre^{-t}, 
    ~~i \,=\, 1,\ldots,n,
    \]
we obtain the $n$th diagonal element of the matrix $\mre ^{-L t} \mre^{-L^T t}$:
    \[
    \left(
    \mre ^{-L t} \mre^{-L^T t}
    \right)_{n}
    \,=\,
    \sum_{i=1}^{n} 
    \left( \frac{t^{i-1}}{(i-1)!} \mre^{-t} \right)^2.
    \]    
Performing the integration from $0$ to $\infty$ yields
\begin{subequations}
  \nonumber
  \begin{align}
    P_{n}
    \,&=\,
    \sum_{i=1}^n \frac{1}{((i-1)!)^2}
   \int_0^\infty
   \frac{\tau^{2(i-1)} \mre^{-\tau}}{2^{2i-1}} 
   \mrd \tau \\    
    \,&=\,
    \sum_{i=1}^n \frac{1}{((i-1)!)^2} \cdot
   \frac{\Gamma(2i-1)}{2^{2i-1}},
  \end{align}
\end{subequations}
where we have used the change of variable $\tau = 2 t$ and the formula for the Gamma function
\begin{equation}
  \label{eq.gamma}
  \Gamma(z) \,=\, \int_0^\infty t^{z-1} \mre^{-\tau} \mrd \tau.
\end{equation}
Since $\Gamma(z) = (z-1)!$, we have the desired formula~\eqref{eq.pn}.

To show the asymptotic scaling of $P_n$, we use Stirling's formula 
\begin{equation}
  \label{eq.stirling}
    n! 
    \,\approx\, 
    \sqrt{2\pi n} 
    \left( \frac{n}{\mre} \right)^n.
\end{equation}
With some algebra, we get 
$
    P_{n}
    \, \approx \,
    \sqrt{n/\pi}.
$
Summing $P_n$ with respect to $n$ yields the expression for the average variance
    $
        \Pi_N
        \,=\,
        \frac{(2N+1)!}{3 \cdot 2^{2N}  N! N!}.
    $
By applying Stirling's formula, we obtain
$    \Pi_N
    \, \approx \,
    \frac{2}{3} \sqrt{N/\pi}.
$

\subsection{Proof of Lemma~\ref{lem.2d}}
\label{app.2d}
Since $L_2$ is lower triangular block Toeplitz, it follows that $(sI + L_2)^{-1}$ is also lower triangular block Toeplitz. In particular,  $(s I + K_2)^{-i}$ is the $i$th block entry of the first column. By using the inverse Laplace transform
    \[
        {\cal L}^{-1}\{
        (sI + K_2)^{-i}\}
        \, = \,
        \frac{t^{i-1} \, \mre^{- K_2 t}}{(i-1)!},
        \quad
        i \, = \, 1,\ldots,N,
    \]
we obtain the $n$th diagonal block of $\mre ^{- L t} \mre^{- L^T t}$, that is,
    \[
        \left(
        \mre ^{- L t} \, \mre^{- L^T t}
        \right)_{n}
        \, = \,
        \sum_{i = 1}^{n}
        \frac{t^{2i-2} }{((i-1)!)^2} \, \mre^{- K_2 t} \, \mre^{- K_2^T t}.
    \]
An analogous calculation shows that the $m$th diagonal element of $\mre^{- K_2 t} \, \mre^{- K_2^T t}$ is given by
    \[
        \left(
        \mre ^{- K_2 t} \, \mre^{- K_2^T t}
        \right)_{m}
        \, = \,
        \sum_{j = 1}^{m}
        \frac{t^{2j-2} }{((j-1)!)^2} \, \mre^{- 2 t} \, \mre^{- 2 t}.
    \]
It follows that
\[
  (P_n)_m \,=\, 
  \int_0^\infty 
  \sum_{i=1}^n \sum_{j=1}^m
  \frac{t^{2i-2}}{((i-1)!)^2}
  \frac{t^{2j-2}}{((j-1)!)^2} \, \mre^{-4 t}  \,
  \mrd t.  
\]
Performing the integration yields the desired formula~\eqref{eq.pnm}.

\subsection{Proof of Lemma~\ref{lem.3d}}
\label{app.3d}
We begin with the steady-state covariance matrix 
    \[
        P
        \, = \,
        \int_{0}^\infty        
        \mre ^{- L_3 t} \, \mre^{- L_3^T t}
        \, \mrd t.
    \]
The matrix exponential 
$
  \mre^{-L_3 t} \,=\, {\cal L}^{-1}\{sI + L_3 \} \in \R^{N^3 \times N^3}
$
is lower triangular block Toeplitz with the $i$th block of the first column being
\[
  {\cal L}^{-1} \{ (sI + K)^{-i} \} 
  \,=\, 
  \frac{t^{i-1}}{(i-1)!} \mre^{-Kt} \in \R^{N^2 \times N^2}.
\]
Since the $j$th block of the first column of $\mre^{-K t}$ is
\[
   {\cal L}^{-1} \{ (sI + K_3)^{-j} \} \,=\, \frac{t^{j-1}}{(j-1)!} \mre^{-K_3 t} 
   \in \R^{N \times N},
\]
and since the $k$th element of the first column of $\mre^{-K_3 t}$ is 
\[
   {\cal L}^{-1} \{ (s + 3)^{-k} \} \,=\, \frac{t^{k-1}}{(k-1)!} \mre^{-3 t},
\]
it follows that
\[
  ((P_n)_m)_l 
  = \int_0^\infty \sum_{i=1}^n   \sum_{j=1}^m   \sum_{k=1}^l 
  \frac{t^{2i+2j+2k-6} \cdot \mre^{-6t}}{((i-1)!(j-1)!(k-1)!)^2}
  \mrd t.  
\]
Performing the integration yields the desired formula~\eqref{eq.3d}.

\subsection{Proof of Proposition~\ref{pro.2d}}
\label{app.2dscale}

\begin{figure}
      \centering
        \includegraphics[width=0.25\textwidth]{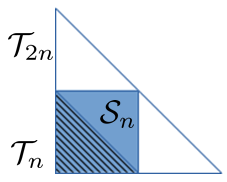}
      \caption{The triangle ${\cal T}_{n}$ (filled in with 45-degree lines) is a subset of the square ${\cal S}_{n}$ (solid square) which itself is a subset of the triangle ${\cal T}_{2n}$ (empty triangle).}
      \label{fig.triangle}
    \end{figure}

We begin by writing $V_{n}$ as 
\[
  V_{n}
  \, = \,
  \sum_{i =0}^{n-1}
  \sum_{j =0}^{n-1}
  f(i,j),
\]
where 
\[
  f(i,j)
  \, = \,
  \frac{(2i+2j)!}{4 \cdot 4^{2(i+j)} \, i! \, i! \, j! \, j!}.
\]
In other words, $V_n$ is the summation of a positive function $f$ over the square ${\cal S}_{n} \DefinedAs \{(i,j) ~|~ 0 \leq i,j \leq n-1\}$. Let $\Delta_n$ be the summation of $f$ over the triangle 
\begin{equation}
  \label{eq.triangle}
  {\cal T}_{n} \,\DefinedAs \, \{(i,j) ~|~ 0 \leq i \leq n-1, 0 \leq j \leq n-i \}, 
\end{equation}
with vertices $(0,0)$, $(0,n-1)$, and $(n-1,0)$:
    \[
        \Delta_n \, \DefinedAs \,
        \sum_{i  =  0}^{n-1}
        \sum_{j  =  0}^{n-i}
        f(i,j).
    \]
Then
    $
        \Delta_{n}
        \, < \,
        V_{n}
        \, < \,
        \Delta_{2n},
    $
because the triangle ${\cal T}_{n}$ is a subset of the square ${\cal S}_{n}$ which itself is a subset of the triangle ${\cal T}_{2n}$; see Fig.~\ref{fig.triangle} for an illustration.

To show 
    $
        V_n \, \sim \, O(\log(n))
    $
for large $n$, it suffices to show
    $
        \Delta_n \, \sim \, O(\log(n)).
    $
We compute the summation of $f$ along the line segment $i+j=k$:
    \begin{align}
    \nonumber
        S_k
       \, & \DefinedAs \, \sum_{i=0}^k f(i,k-i)
        \, = \,
        \sum_{i = 0}^k \frac{(2k)!}{4 \cdot 4^{2k} \; i!\,i!\,(k-i)!\,(k-i)!}
        \\
    \nonumber
        \, &= \,
        \frac{1}{4 \cdot 4^{2k}}
        \frac{(2k)!}{k!k!}  \sum_{i=0}^k
        \frac{k!k!}{i!\,i!\,(k-i)!\,(k-i)!} \\
    \label{eq.Sk}
        \, &= \,
        \frac{1}{4 \cdot 4^{2k}} \left( \frac{(2k)!}{k!k!} \right)^2,
        \quad
        k=0,1,\ldots,n-1,
    \end{align}
where we have used the fact that $\sum_{i=0}^k {k \choose i}^2 = {2k \choose k} = \frac{(2k)!}{k!k!}$. From the expression~\eqref{eq.pn} and the approximation $P_n \approx \sqrt{n/\pi}$ in Lemma~\ref{lem.1d}, we conclude that for large $k$,
\[
  S_k 
  \, \approx \, \frac{1}{4\pi k}.
\]
It follows that
    \[
        \Delta_n \,=\, \sum_{k=0}^{n-1} S_k
        \, \sim \, O(\log(n)).
    \]
This completes the proof.

\subsection{Proof of Proposition~\ref{pro.3d}}
\label{app.3dscale}

Setting $n=m=l$ in~\eqref{eq.3d} yields the variance of the follower $(n,n,n)$:
\[
   V_n \,=\, 
    \sum_{i=0}^{n-1}\sum_{j=0}^{n-1}\sum_{k=0}^{n-1} g(i,j,k),
\]
where
\[
  g(i,j,k) \,=\, \frac{(2i+2j+2k)!}{6 \cdot 6^{2(i+j+k)} (i! j! k!)^2}.
\]
In other words, $V_n$ is the summation of the positive function $g$ over the cube ${\cal C}_n \DefinedAs \{0 \leq i,j,k \leq n-1\}$. Let $T_n$ be the summation of $g$ over the triangular pyramid ${\cal P}_n \DefinedAs \{0 \leq p \leq n-1, 0\leq j \leq p, 0 \leq k \leq p-j\}$, whose vertices are given by $\{(0,0,0), (n-1,0,0), (0,n-1,0), (0,0,n-1)\}$. It follows that
$
  T_n \,<\, V_n \,<\, T_{2n}.
$
This is because ${\cal P}_n$ is a subset of ${\cal C}_n$ which itself is a subset of ${\cal P}_{2n}$. Thus, it suffices to show that $T_n \sim O(1)$. 

We compute the sum of $g$ across the triangle segment of the pyramid
$
  T_n \,=\, \sum_{p=0}^{n-1} G_p,
$
where $G_p$ is the sum of $g$ over the triangle $i+j+k=p$:
  \begin{align}
     \nonumber
     G_p &=\, 
         \sum_{j=0}^p \sum_{k=0}^{p-j} f(p-j-k,j,k) \\
    \label{eq.deltap}
         &=\,
         \sum_{j=0}^p \sum_{k=0}^{p-j} 
         \frac{(2p)!}{6\cdot 2^{2p} p!p!} 
         \left(\frac{p!}{3^{p} j! k! (p-j-k)!} \right)^2.
  \end{align}
To evaluate the summation, we employ a probability argument from~\cite{doysne84}. Consider dropping $p$ balls into three boxes $A$, $B$, and $C$. The probability of dropping $j$ balls into $A$, $k$ balls into $B$, and $p-j-k$ balls into $C$ is $\frac{p!}{3^{p} j! k! (p-j-k)!}$. Since the largest probability occurs when the same number of balls drop in three boxes, it follows that
\[
  G_p  \leq 
         \frac{(2p)!}{6\cdot 2^{2p} p!p!} \cdot
    \frac{p!}{3^p (\lfloor \frac{p}{3} \rfloor!)^3} 
    \sum_{j=0}^p \sum_{k=0}^{p-j} \left( \frac{p!}{3^{p} j! k! (p-j-k)!} \right)^2,
\]
where $ \lfloor \frac{p}{3} \rfloor$ denotes the largest integer that is no greater than $p/3$. Note that 
\[
  \sum_{j=0}^p \sum_{k=0}^{p-j} \left( \frac{p!}{3^{p} j! k! (p-j-k)!} \right)^2 \,=\, 1,
\]
because it is the sum of the probability of all outcomes of dropping three balls in three boxes. Therefore, for large $n$,
\[
  G_p \, \leq \, 
         \frac{1}{6\cdot 2^{2p}}  \cdot
         \frac{(2p)!}{p!p!}  \cdot
         \frac{p!}{3^p (\lfloor \frac{p}{3} \rfloor!)^3}  
  \, \approx \, 
  c \, p^{-3/2},
\]
where $c$ is a constant and  we have used Stirling's formula~\eqref{eq.stirling}. It follows that 
\[
  T_n 
  \, = \, \sum_{p=1}^n G_p 
  \, \approx \, \sum_{p=1}^n c \, p^{-3/2} 
  \, \sim \, O(1).
\]
This completes the proof.

\section{Conclusions}
\label{sec.concl}

In this paper, we have obtained explicit formulas for the steady-state variance distribution of leader-follower multi-agent systems in directed lattices in one, two, and three dimensions. We show that the variance of the followers scales as a square-root function of the distance from the leader in the 1D lattice, scales as a logarithmic function along the diagonal of the 2D lattice, and is bounded by a network-size independent constant in the 3D lattice. 

\bibliographystyle{IEEEtran}
\bibliography{network_resistance}

\end{document}